\documentstyle[12pt,us,titlepage]{article}

\begin{document}
\def\l{\lambda}
\def\m{\mu}
\def\b{\beta}
\def\d{\Delta}
\begin{titlepage}
\vspace*{\fill}
\begin{center}
\large\bf NUMERICAL COMPUTATION OF REAL OR COMPLEX ELLIPTIC INTEGRALS
\end{center}
\vfill
\begin{center}
{\large B.~C.~Carlson}
\vfill
{\normalsize
Ames Laboratory and Department of Mathematics, Iowa State
 University,\\
 Ames, Iowa 50011-3020, USA}
\end{center}
\vfill
{\bf Abstract.}\quad Algorithms for numerical computation of symmetric 
elliptic integrals of all three kinds are improved in several ways and 
extended to complex values of the variables (with some restrictions in the
case of the integral of the third kind).  Numerical check values, consistency 
checks, and relations to Legendre's integrals and Bulirsch's integrals are
included.
\vfill
\\
\\
{\bf AMS(MOS) subject classifications.}\quad primary 33A25, 33-04, 65D20;
secondary 33A10, 30-04, 32-04 
\\
\\
This work was supported by the Director of Energy Research, Office of Basic
Energy Sciences.  The Ames Laboratory is operated for the U.~S.~Department of
Energy by Iowa State University under Contract W-7405-ENG-82.

\end{titlepage}
\vspace*{\fill}
\begin{center}
Running \vspace{1.0in} title\\
{\sc\bf COMPUTATION OF ELLIPTIC INTEGRALS}
\end{center}
\vfill
\newpage
\section{Introduction}
\indent Let \(f(x)\) be a real function that is rational except for the square 
root of a cubic or quartic polynomial with at least one pair of conjugate 
complex zeros.  Then \(\int f(x) dx\) can be expressed in terms of standard 
elliptic integrals with complex variables, which are subsequently changed to 
real variables by using quadratic transformations \cite{Ca91}\cite{Ca92}.
Since the transformations complicate the formulas, it is desirable to have 
algorithms for numerical computation of standard elliptic integrals with 
complex variables.  Such integrals are met in other contexts also, for example 
in conformal mapping, and the complex variables might not occur in conjugate 
pairs.

This paper contains algorithms for numerical computation of complete and 
incomplete elliptic integrals of all three kinds when the variables are complex
(with some restrictions for integrals of the third kind).  They are similar
to algorithms published earlier \cite{Ca79} for real variables,
but several improvements that apply to the real case have been made in the 
course of extending them to complex variables.  The integrals computed are
the symmetric integrals of the first and third 
kinds and two degenerate cases, of which one is an elementary function and the 
other is an elliptic integral of the second kind.  Other integrals can be 
obtained from these by using the formulas and references in Section 4.  The 
method of computation is to iterate the duplication theorem and then sum a 
power series up to terms of degree five; the error ultimately decreases by 
a factor of \(4^6 = 4096\) with
each duplication.  Since this method is slowest when the integrals are 
complete, we add a faster algorithm for computing complete integrals of the
first and second kinds, including Legendre's \(K(k)\) and \(E(k)\) for complex
\(k\), by the method of arithmetic and geometric means.

The symmetric integral of the first kind is
\begin{equation}
\label{1.1}
R_F(x,y,z) = \frac{1}{2}\int_{0}^{\infty}[(t+x)(t+y)(t+z)]^{-1/2}dt\,,
\end{equation}
where the square root is taken real and positive if \(x,y,z\) are positive
and varies continuously when \(x,y,z\) become complex.  The integral
is well defined if \(x,y,z\) lie in the complex plane cut along the 
nonpositive real axis (henceforth called the ``cut plane''), with the 
exception that at most one of \(x,y,z\) may be 0.
The same requirements apply to the symmetric integral of the third kind,
\begin{equation}\label{1.2}
R_J(x,y,z,p) =
\frac{3}{2}\int_{0}^{\infty}[(t+x)(t+y)(t+z)]^{-1/2}(t+p)^{-1}dt\,,
\end{equation}
where \(p \neq 0\) and the Cauchy principal value is to be taken if \(p\) is 
real and negative.  A degenerate case of \(R_F\) that embraces the inverse 
circular and inverse hyperbolic functions (see Section 4) is 
\begin{equation}\label{1.3}
R_C(x,y) = R_F(x,y,y) =
\frac{1}{2}\int_{0}^{\infty}(t+x)^{-1/2}(t+y)^{-1}dt\,.
\end{equation}
It is well defined if \(x\) lies in the cut plane or is 0 and if \(y \neq 0\);
the Cauchy principal value is to be taken if \(y\) is real and negative.
Professor Luigi Gatteschi pointed out to me that Fubini, while still a student
at Pisa, proposed the use of a function equivalent to \(1/R_C\)
in his first published paper \cite{Fu97}.  A degenerate case of \(R_J\) that is
an elliptic integral of the second kind is
\begin{equation}\label{1.4}
R_D(x,y,z) = R_J(x,y,z,z) =
\frac{3}{2}\int_{0}^{\infty}[(t+x)(t+y)]^{-1/2}(t+z)^{-3/2}dt\,,
\end{equation}
which is well defined under the same conditions as \(R_F\) except that \(z\) 
must not be 0.

     Because \(R_D\) is symmetric only in \(x\) and \(y\), it is sometimes 
convenient to use a completely symmetric integral of the second kind,
\begin{equation}\label{1.5}
R_G(x,y,z) = \frac{1}{4}\int_{0}^{\infty}[(t+x)(t+y)(t+z)]^{-1/2}
\left( \frac{x}{t+x} + \frac{y}{t+y} + \frac{z}{t+z} \right) tdt\,,
\end{equation}
where any or all of \(x,y,z\) may be 0 and those that are nonzero lie in the 
cut plane.  If the closed convex hull of 
\( \{x,y,z\} \) lies in the union of 0 and the cut plane, \(R_G\) is 
represented by a double integral that accounts for its usefulness in 
problems connected with ellipsoids,
\begin{equation}\label{1.6}
   R_G(x,y,z) = \frac{1}{4\pi} \int_0^{2\pi} \int_0^{\pi} (x \sin^2\!\theta
   \cos^2\!\phi + y \sin^2\!\theta \sin^2\!\phi + z \cos^2\!\theta)^{1/2}
   \sin\theta\,d\theta\,d\phi\,.
\end{equation}
Except that at most one of \(x,y,z\) may be 0, \(R_F\) has the same 
representation with exponent \(-1/2\) instead of 1/2.  With the same 
exception and with \(x,y,z\) permuted so that \(z \neq 0\),
\(R_G\) can be obtained from \(R_F\) and \(R_D\) by using the relation 
\begin{equation}\label{1.7}
   2R_G(x,y,z) = zR_F(x,y,z) - \frac13 (x-z)(y-z)R_D(x,y,z)
   + \sqrt{\frac{xy}{z}}.
\end{equation}
Algorithms for direct computation of \(R_G\) and \(R_F\) with real variables
by successive Landen transformations are given in \cite{Ca65}.  In Section 2 of
the present paper, one of those algorithms is extended to complex variables, 
but only in the complete case.  For each of the functions defined above, the 
complete case is the case in which one of \(x,y,z\) is 0.

The functions \(R_F\) and \(R_C\) are homogeneous of degree \(-1/2\) in their 
variables, \(R_J\) and \(R_D\) of degree \(-3/2\), and \(R_G\) of degree +1/2.
Their homogeneity and symmetry replace a set of linear transformations of 
Legendre's integrals, simplify their quadratic transformations and other 
properties, and make it possible (see \cite{Ca87}-\cite{Ca92})
to unify many of the formulas in customary integral tables.

The earlier versions of the algorithms in Section 2 (except the last one)
were published in 
\cite{Ca79}, modified in \cite[(4),(5)]{CN81} to avoid underflows, and 
implemented by Fortran codes in several major software libraries, in the 
Supplements to \cite{Ca87} and \cite{Ca88}, and in \cite{PT90} and 
\cite[\S6.11]{PT92a}.  Codes in C can be found in \cite[\S6.11]{PT92b}.  The 
algorithms in the present paper are preferable to those in \cite{Ca79} in the 
following four respects (in addition to incorporating the modification in 
\cite{CN81}):

(1) All complex values of the variables are admissible for which \(R_F, R_D\,
\), and \(R_C\) were defined above, with the exception that computing the 
Cauchy principal value of \(R_C(x,y)\) when \(y < 0\) requires the preliminary 
transformation (\ref{2.13a}).  The variables of \(R_J\) are 
restricted by conditions that are sufficient but not necessary to keep the 
fourth variable from being transformed to 0 by the duplication theorem.  

(2) The algorithms have been rearranged to reduce substantially the number of
arithmetic operations.

(3) A bound on the fractional error of the result can be specified directly.

(4) Because \(R_C\) is used repeatedly in the algorithm for \(R_J\), its 
computation has been speeded up a little by including terms up to degree seven 
(instead of five) in the truncated power series.

Section 3 contains information for checking codes based on the algorithms of 
Section 2, and Section 4 relates other integrals to the integrals used here.

\section{Algorithms}
We shall summarize briefly the method of computation; details are given
in \cite{Ca79} for the case of real variables, and a full discussion of the 
complex case with proof of error bounds will appear eventually in a book now 
in preparation.  Validity of algorithms can be tested in the complex domain
by the consistency checks in Section 3 and sometimes by comparison of the 
complete case with the last algorithm in the present Section.

If we obtain \(x_{m+1}\,,\,y_{m+1}\,,\,z_{m+1}\) from \(x_m\,,\,y_m\,,\,z_m\) 
as prescribed by (\ref{2.4}) below, the duplication theorem implies
\begin{equation}\label{2.1}
   R_F(x_0\,,\,y_0\,,\,z_0) = R_F(x_m\,,\,y_m\,,\,z_m)\,,\qquad m = 
   1,2,3,\ldots\,.
\end{equation}
Although \(x_m\,,\,y_m\,,\,z_m\) converge to a common limit in the cut plane
as \(m \rightarrow \infty\), differences like \(x_m - y_m\) decrease by a 
factor of only 4 when \(m\) increases by 1.  When the fractional differences 
between \(x_m\,,\,y_m\,,\,z_m\) and their arithmetic average become smaller 
than an amount determined by the desired accuracy of the final result, \(R_F\) 
is expanded in a multiple series of powers of these fractional differences, 
denoted by \(X,Y,Z\).  Because of the symmetry of \(R_F\) the series can be 
rewritten in terms of the elementary symmetric functions 
\(E_1\,,\,E_2\,,\,E_3\) of \(X,Y,Z\); and because \(E_1 = X+Y+Z = 0\) the terms
up to degree 5 are very simple.  The truncation error, being of degree 6, 
ultimately decreases by a factor of \(4^6\) with each duplication.  An estimate
of truncation error is provided by \cite[(A.10)]{Ca79}.
\\
\\
\noindent {\sc ALGORITHM FOR $R_F\,$}.
\begin{em}
We suppose that at most one of \(x,y,z\) is 0 and those that are nonzero have
phase less in magnitude than \(\pi\).
The function \(R_F(x,y,z)\) defined by (\ref{1.1}) is to be computed with
relative error less in magnitude than r.  (We assume \(r < 3 \times 10^{-4}\).)

Let \(x_0 = x,\, y_0 = y,\, z_0 = z,\,\) and
\begin{equation}\label{2.2}
     A_0 = \frac{x+y+z}{3}\,,\qquad Q = (3r)^{-1/6}\max\{|A_0 - x|,|A_0 - y|,
     |A_0 - z|\}\,.
\end{equation}
For \(\ m = 0,1,2,\ldots\,\), define
\begin{equation}\label{2.3}
   \l_m = \sqrt{x_m}\,\sqrt{y_m} + \sqrt{x_m}\,\sqrt{z_m} +
    \sqrt{y_m}\,\sqrt{z_m}\,,\qquad A_{m+1} = \frac{A_m + \l_m}{4}\,,
\end{equation}
\begin{equation}\label{2.4}
   x_{m+1}=\frac{x_m + \l_m}{4}\,,\qquad y_{m+1}=\frac{y_m + \l_m}{4}\,,\qquad
   z_{m+1} = \frac{z_m + \l_m}{4}\,,
\end{equation}
where each square root has nonnegative real part.

Compute \(A_m\;\)for \(\ m = 0,1,\ldots,n\,\), where
\( 4^{-n}Q < |A_n|\,\).  Define
\begin{equation}\label{2.5}
   X = \frac{A_0 - x}{4^n A_n}\,,\qquad Y = \frac{A_0 - y}{4^n A_n}\,,\qquad
   Z = -X - Y\,,
\end{equation}
\begin{equation}\label{2.6}
   E_2 = XY - Z^2\,,\qquad E_3 = XYZ\,.
\end{equation}
Then
\begin{equation}\label{2.7}
   R_F(x,y,z) \approx A_n^{-1/2}\left( 1 - \frac{1}{10}E_2 + \frac{1}{14}E_3
   + \frac{1}{24}E_2^2 - \frac{3}{44}E_2E_3 \right)
\end{equation}
with relative error less in magnitude than r.
\end{em}
\\
\\
\indent The statement that the relative error is less in magnitude than \(r\) 
means that the true value of the function lies inside a circle in the complex 
plane with center at the computed value and radius equal to \(r\) times the
distance of the center from the origin.  We assume that \(r\) is large compared
to the machine precision, so that roundoff error is negligible compared to the 
error produced by approximations made in the algorithm.  These remarks apply to
all the algorithms in this paper.

We note the relations
\begin{equation}\label{2.8}
   A_m = \frac{x_m + y_m + z_m}{3}\,,\qquad A_m - x_m =
   \frac{A_0 - x}{4^m}\,,\qquad X = 1 - \frac{x_n}{A_n}\,,
\end{equation}
and similar relations obtained by permuting \((x,X),(y,Y),(z,Z)\,\).
Although \(A_0\) can be 0 in the complex case, \(A_n \neq 0\) because
\(Q \geq 0\).  (It can be shown that \(A_m \neq 0\) if \(m \geq 1\).)
If \(x_m\,,\,y_m\,,\,z_m\,\) are real and nonnegative, the inequality of
arithmetic and geometric means implies \(\l_m \leq 3A_m\,\), whence
\(A_{m+1} \leq A_m\,\).

If \(y = z\) then \(R_F\) reduces to \(R_C\).  Because \(E_2\) and \(E_3\) are 
no longer independent, the series in (\ref{2.7}) becomes a series in one 
variable.  By including terms up to degree seven (instead of five), we can 
usually save one duplication, which seems worthwhile because computation of 
\(R_J\) requires one computation of \(R_C\) in each cycle of iteration.
\\
\\
\noindent {\sc ALGORITHM FOR $R_C\,$}.
\begin{em}
Let \(x\) and \(y\) be nonzero and have phase less in magnitude
than \(\pi\), with the exception that \(x\) may be 0.
The function \(R_C(x,y)\) defined by (\ref{1.3}) is to be computed
with relative error less in magnitude than \(r\). (We assume \(r < 2 \times
10^{-4}\).)

Let \(x_0 = x\,,\,y_0 = y\,\), and
\begin{equation}\label{2.9}
   A_0 = \frac{x+2y}{3}\,,\qquad\quad Q = (3r)^{-1/8} |A_0 - x|\,.
\end{equation}
For \(\;m = 0,1,2,\ldots\,,\) define
\begin{equation}\label{2.10}
   \l_m = 2\sqrt{x_m}\,\sqrt{y_m} + y_m\,,\qquad\quad A_{m+1} =
   \frac{A_m + \l_m}{4}\,,
\end{equation}
\begin{equation}\label{2.11}
x_{m+1} = \frac{x_m + \l_m}{4}\,,\qquad\quad y_{m+1} = \frac{y_m + \l_m}{4}\,,
\end{equation}
where each square root has nonnegative real part.

Compute \(A_m\;\) for \(\ m = 0,1,\ldots,n\,\), where \( 4^{-n}Q < |A_n|\,\),
and define
\begin{equation}\label{2.12}
   s = \frac{y - A_0}{4^n A_n}\,.
\end{equation}
Then
\begin{equation}\label{2.13}
   R_C(x,y) \approx A_n^{-1/2} \left( 1 +
   \frac{3}{10}s^2 + \frac{1}{7}s^3 + \frac{3}{8}s^4 + \frac{9}{22}s^5
   + \frac{159}{208}s^6 + \frac{9}{8}s^7 \right)
\end{equation}
with relative error less in magnitude than \(r\).
\end{em}
\\
\\
\indent If the second variable of \(R_C\) is real and negative, the Cauchy 
principal value is
\begin{equation}\label{2.13a}
R_C(x,-y) = \left(\frac{x}{x+y}\right)^{1/2} R_C(x+y,y)\,, \qquad y>0\,,
\end{equation}
by \cite[(4.8)]{ZC} and (\ref{4.20}).  This vanishes if \(x=0\).

In the notation of (\ref{2.17a}) to (\ref{2.19}), the duplication 
theorem for \(R_J\) is
\begin{equation}\label{2.14}
   R_J(x_m,y_m,z_m,p_m) = \frac{1}{4} R_J(x_{m+1},y_{m+1},z_{m+1},p_{m+1})
   + \frac{6}{d_m}\,R_C(1,\,1+e_m)\,.
\end{equation}
The duplication theorem for \(R_C\) has been applied to \cite[(5.1)]{Ca79} 
to allow a wider range of complex phase for the variables.  Iteration of
(\ref{2.14}) yields
\begin{equation}\label{2.15}
   R_J(x_0\,,\,y_0\,,\,z_0\,,\,p_0) = 4^{-n}\,R_J(x_n,y_n,z_n,p_n) + 
   6 \sum_{m=0}^{n-1} \frac{4^{-m}}{d_m}\,R_C(1,\,1+e_m)\,.
\end{equation}
The first term on the right side can be treated as a symmetric function of the 
five variables \(x_n\,,\,y_n\,,\,z_n\,,\,p_n\,,\,p_n\) by writing 
\((t+p)^{-1}\) in (\ref{1.2}) as \((t+p)^{-1/2} (t+p)^{-1/2}\).  When the 
fractional differences \(X,Y,Z,P,P\) between these five variables and 
their arithmetic average become small
enough, \(R_J\) is expanded in powers of the elementary symmetric functions 
\(E_2\,,\,E_3\,,\,E_4\,,\,E_5\) of \(X,Y,Z,P,P\) (since \(E_1 = 0\) ).
\\
\\
\noindent {\sc ALGORITHM FOR $R_J\,$}.
\begin{em}
Let \(x,y,z\) have nonnegative real part and at most one of them be
0, while {\rm Re}\(\,p > 0\).  Alternatively, if \(p \neq 0\) and 
\( |{\rm ph}\,p| < \pi\), either let \(x,y,z\) be real and nonnegative and 
at most one of them be 0, or else let two of the variables \(x,y,z\) be 
nonzero and conjugate complex with phase less in magnitude than \(\pi\) 
and the third variable be real and nonnegative.  The function
\(R_J(x,y,z,p)\) defined by (\ref{1.2}) is to be computed with relative error
less in magnitude than \(r\).  (We assume \(r < 10^{-4}\).)

Let \((x_0,y_0,z_0,p_0) = (x,y,z,p)\) and
\begin{equation}\label{2.16}
     A_0 = \frac{x+y+z+2p}{5}\,,\qquad \delta = (p-x)(p-y)(p-z)\,,
\end{equation}
\begin{equation}\label{2.17}
    Q = (r/4)^{-1/6}\max\{|A_0 - x|,|A_0 - y|,|A_0 - z|,|A_0 - p|\}\,.
\end{equation}
For \(\ m = 0,1,2,\ldots\,\), define
\begin{equation}\label{2.17a}
   \l_m = \sqrt{x_m}\,\sqrt{y_m} + \sqrt{x_m}\,\sqrt{z_m} +
   \sqrt{y_m}\,\sqrt{z_m}\,,\qquad A_{m+1} = \frac{A_m + \l_m}{4}\,,
\end{equation}
\begin{equation}\label{2.18}
   x_{m+1}=\frac{x_m + \l_m}{4}\,,\quad y_{m+1}=\frac{y_m + \l_m}{4}\,,\quad
   z_{m+1} = \frac{z_m + \l_m}{4}\,,\quad p_{m+1} = \frac{p_m + \l_m}{4}\,,
\end{equation}
\begin{equation}\label{2.19}
   d_m = (\sqrt{p_m}\,+\sqrt{x_m}\,)(\sqrt{p_m}\,+\sqrt{y_m}\,)
   (\sqrt{p_m}\,+\sqrt{z_m}\,)\,,\qquad e_m = \frac{4^{-3m}\,\delta}{d_m^2}\,,
\end{equation}
where each square root has nonnegative real part.

Compute \(A_m\;\)for \(\ m = 0,1,\ldots,n\,\), where \(4^{-n}Q < |A_n|\).
Compute also \(R_C(1,\,1+e_m)\,\) with relative error less in magnitude 
than r for \(m = 0,1,\ldots,n-1\).  Define
\begin{equation}\label{2.20}
   X = \frac{A_0-x}{4^n A_n}\,,\qquad Y = \frac{A_0-y}{4^n A_n}\,,\qquad
   Z = \frac{A_0-z}{4^n A_n}\,,\qquad P = (-X-Y-Z)/2\,,
\end{equation}
\begin{eqnarray}
   & & E_2 = XY+XZ+YZ-3P^2,\qquad E_3 = XYZ + 2E_2P + 4P^3,\label{2.21}\\
   & & E_4 = (2XYZ + E_2P + 3P^3)P\,,\qquad E_5 = XYZP^2.\label{2.22}
\end{eqnarray}
Then
\begin{eqnarray}
   && R_J(x,y,z,p) \approx 4^{-n}A_n^{-3/2} \left( 1 - \frac{3}{14}E_2 +
   \frac{1}{6}E_3 + \frac{9}{88}E_2^2 - \frac{3}{22}E_4 - \frac{9}{52}E_2E_3
   + \frac{3}{26}E_5 \right) \nonumber\\
   &&\qquad\qquad\qquad +\; 6 \sum_{m=0}^{n-1} \frac{4^{-m}}{d_m}\,
   R_C(1,\,1+e_m) \label{2.23}
\end{eqnarray}
with relative error less in magnitude than \(r\).
\end{em}
\\
\\
\indent When the variables are complex, the bound on relative error is not 
rigorous because of the possibility of some cancellation between terms 
that individually have error less than \(r\).  In practice, 
however, the error is usually much smaller than the bound.

If \(x,y,z\) are real and nonnegative, at most one of them is 0, and the 
fourth variable of \(R_J\) is negative, 
the Cauchy principal value is given by \cite[(4.6)]{ZC} and (\ref{2.13a})\,:
\begin{eqnarray}
(y+q)R_J(x,y,z,-q) & = & (p-y)R_J(x,y,z,p) - 3R_F(x,y,z) \nonumber\\
& + & 3\left(\frac{xyz}{xz+pq}\right)^{1/2}R_C(xz+pq,pq), \quad q>0,
\label{2.23a}
\end{eqnarray}
where \(p-y = (z-y)(y-x)/(y+q)\).  If we permute the values of \(x,y,z\) so
that \( (z-y)(y-x) \geq 0 \), then \(p \geq y > 0 \) and all terms on the right
side can be computed by the preceding algorithms.  If \(x,y,z\) are complex,
conditions of validity have not been established for this method of computing 
the Cauchy principal value.

It would be desirable to give the algorithm less restrictive conditions 
on \(x,y,z,p\) that do not rule out so many cases in which \(p = z\) 
and \(R_J\) reduces to \(R_D\), for in such cases the algorithm
gives correct results under the much weaker conditions stated below in the
algorithm for \(R_D\).  We note that \(p = z\) implies \(e_m = 0\) and 
\(d_m = 2\sqrt{z_m}(z_m + \l_m)\) for all \(m\), whence (\ref{2.15}) becomes
\begin{equation}\label{2.24}
   R_D(x_0\,,\,y_0\,,\,z_0) = 4^{-n}R_D(x_n,y_n,z_n) 
   + 3 \sum_{m=0}^{n-1}\frac{4^{-m}}{\sqrt{z_m}\,(z_m + \l_m)}\,.
\end{equation}
The function \(R_D(x,y,z)\) is treated as a symmetric function of
\(x,y,z,z,z\).
\\
\\
\noindent {\sc ALGORITHM FOR $R_D\,$}.
\begin{em}
Let \(x,y,z\) be nonzero and have phase less in magnitude than
\(\pi\), with the exception that at most one of \(x\) and \(y\) may be 0.
The function \(R_D(x,y,z)\) defined by (\ref{1.4}) is to be computed with
relative error less in magnitude than r.  (We assume \(r < 10^{-4}\).)

Let \(x_0 = x,\, y_0 = y,\, z_0 = z,\,\) and
\begin{equation}\label{2.25}
     A_0 = \frac{x+y+3z}{5}\,,\qquad Q = (r/4)^{-1/6}\max\{|A_0 - x|,|A_0 - y|,
     |A_0 - z|\}\,.
\end{equation}
For \(\ m = 0,1,2,\ldots\,\), define
\begin{equation}\label{2.26}
   \l_m = \sqrt{x_m}\,\sqrt{y_m} + \sqrt{x_m}\,\sqrt{z_m} +
   \sqrt{y_m}\,\sqrt{z_m}\,,\qquad A_{m+1} = \frac{A_m + \l_m}{4}\,,
\end{equation}
\begin{equation}\label{2.27}
   x_{m+1}=\frac{x_m + \l_m}{4}\,,\qquad y_{m+1}=\frac{y_m + \l_m}{4}\,,\qquad
   z_{m+1} = \frac{z_m + \l_m}{4}\,,
\end{equation}
where each square root has nonnegative real part.

Compute \(A_m\;\)for \(\ m = 0,1,\ldots,n\,\), where
\( 4^{-n}Q < |A_n|\,\).  Define
\begin{equation}\label{2.28}
   X = \frac{A_0 - x}{4^n A_n}\,,\qquad Y = \frac{A_0 - y}{4^n A_n}\,,\qquad
   Z = -(X+Y)/3\,,
\end{equation}
\begin{eqnarray}
   && E_2 = XY-6Z^2\,,\quad E_3 = (3XY-8Z^2)Z\,,\label{2.29}\\
   && E_4 = 3(XY-Z^2)Z^2\,,\qquad E_5 = XYZ^3.\label{2.30}
\end{eqnarray}
Then
\begin{eqnarray}
   && R_D(x,y,z) \approx 4^{-n}A_n^{-3/2} \left( 1 - \frac{3}{14}E_2 +
   \frac{1}{6}E_3 + \frac{9}{88}E_2^2 - \frac{3}{22}E_4 - \frac{9}{52}E_2E_3
   + \frac{3}{26}E_5 \right) \nonumber\\
   &&\qquad\qquad\qquad +\; 3 \sum_{m=0}^{n-1} \frac{4^{-m}}{\sqrt{z_m}\,
   (z_m + \l_m)}\label{2.31}
\end{eqnarray}
with relative error less in magnitude than \(r\).
\end{em}
\\
\\
\indent The remark about relative error immediately following the algorithm
for \(R_J\) applies again here.  The function \(R_G\) can be computed from
\(R_F\) and \(R_D\) by using (\ref{1.7}) if at most one of \(x,y,z\) is 0.
Neither (\ref{1.7}) nor the next algorithm can be used 
to compute \(R_G(0,0,z) = \sqrt{z}\,/2\).

A faster and simpler way of computing the complete case of \(R_F\) and \(R_G\) 
with complex variables is useful because Legendre's complete integrals of the 
first and second kinds are
\begin{equation}\label{2.32}
   K(k) = R_F(1-k^2,\,1\,,\,0),\qquad  E(k) = 2R_G(1-k^2,\,1\,,\,0).
\end{equation}
The method of successive arithmetic and geometric means (a special case of
successive Landen transformations) serves this purpose for
complex \(k^2 \notin [1,+\infty) \).  (Algorithms using Landen transformations
for incomplete \(R_F\) and 
\(R_G\) are given in \cite{Ca65} but only for real variables; conditions of 
validity for complex variables are not obvious.)
\\
\\
\noindent {\sc ALGORITHM FOR $R_F(x\,,\,y\,,\,0)$ AND $R_G(x\,,\,y\,,\,0)$}.
\begin{em}
Let \(x\) and \(y\) be nonzero and have phase less in magnitude than \(\pi\).
The complete elliptic integrals \(R_F(x\,,\,y\,,\,0)\)
and \(R_G(x\,,\,y\,,\,0)\) are to be computed with relative error less in
magnitude than \(r\).

Define \(x_0 = \sqrt{x}\,,\,y_0 = \sqrt{y}\,,\) and
\begin{equation}\label{2.33}
   x_{m+1} = \frac{x_m + y_m}{2}\,,\,\quad y_{m+1} = \sqrt{x_m y_m}\,,
   \qquad m = 0,1,2,\ldots\,,
\end{equation}
where each square root has positive real part.  Compute \(x_m\) and \(y_m\) for
m = 0,1,\ldots,n\,, where
\begin{equation}\label{2.34}
   |x_n - y_n| < 2.7\,\sqrt{r}\,|x_n|\,.
\end{equation}
Then
\begin{equation}\label{2.35}
   R_F(x,y,0) \approx \frac{\pi}{x_n + y_n}
\end{equation}
with relative error less in magnitude than r.  Also,
\begin{equation}\label{2.36}
   2R_G(x,y,0) \approx \left( \left( \frac{x_0 + y_0}{2} \right)^2 
   - \sum_{m=1}^n 2^{m-2} (x_m - y_m)^2 \right) R_F(x,y,0)
\end{equation}
with relative error less in magnitude than \(r\) if we neglect terms of order
\(r^2\).  The summation is empty if \(n = 0\).
\end{em}
\\
\\
\indent This algorithm can be used to compute also
\begin{equation}\label{2.37}
   R_D(0,y,z) = \frac{3}{z(y-z)}[2R_G(y,z,0) - z\,R_F(y,z,0)]\,,\qquad 
   0 \neq y\neq z \neq 0\,.
\end{equation}
The exceptional case with \(y=z \neq 0\) is
\begin{equation}\label{2.38}
   R_D(0,\,y,\,y) = \frac{3\pi}{4}\,y^{-3/2}.
\end{equation}
\section{Numerical checks}
\indent Codes based on the algorithms of Section 2 can be checked against the 
following assortment of numerical values for complete and incomplete integrals 
with real, conjugate complex, or nonconjugate complex variables.
\begin{eqnarray*}
   && R_F(1,2,0) = 1.3110\;28777\;1461\\
   && R_F(i,-i,0) = R_F(0.5,1,0) = 1.8540\;74677\;3014\\
   && R_F(i-1,i,0) = 0.79612\;58658\;4234\;-\;i\;1.2138\;56669\;8365\\
   && R_F(2,3,4) = 0.58408\;28416\;7715\\
   && R_F(i,-i,2) = 1.0441\;44565\;4064\\
   && R_F(i-1,i,1-i) = 0.93912\;05021\;8619\;-\;i\;0.53296\;25201\;8635\\
\\
   && R_C(0,1/4) = \pi = 3.1415\;92653\;5898 \\
   && R_C(9/4,2) = \ln 2 = 0.69314\;71805\;5995 \\
   && R_C(0,i) = (1-i)\;1.1107\;20734\;5396 \\
   && R_C(-i,i) = 1.2260\;84956\;9072\;-\;i\;0.34471\;13698\;8768 \\
   && R_C(1/4\,,\,-2) = \frac{\ln 2}{3} = 0.23104\;90601\;8665 \\
   && R_C(i\,,\,-1) = 0.77778\;59692\;0447\;+\;i\;0.19832\;48499\;3429\\
\\
   && R_J(0,1,2,3) = 0.77688\;62377\;8582 \\
   && R_J(2,3,4,5) = 0.14297\;57966\;7157 \\
   && R_J(2,3,4,-1+i) = 0.13613\;94582\;7771\;-\;i\;0.38207\;56162\;4427\\
   && R_J(i,-i,0,2) = 1.6490\;01166\;2711 \\
   && R_J(-1+i,-1-i,1,2) = 0.94148\;35884\;1220 \\
   && R_J(i,-i,0,1-i) = 1.8260\;11522\;9009\;+\;i\;1.2290\;66190\;8643 \\
   && R_J(-1+i,-1-i,1,-3+i) = -0.61127\;97081\;2028\;-\;i\;
      1.0684\;03839\;0007 \\
   && R_J(-1+i,-2-i,-i,-1+i) = 1.8249\;02739\;3704\;-\;i\;1.2218\;47578\;4827
\end{eqnarray*}
\\
The last case does not fit the assumptions, but the algorithm yields a
value agreeing with \(R_D(-2-i,-i,-1+i)\) below.
If \(x,y,z\) are strictly positive, the Cauchy principal value of \(R_J\) 
changes sign at a negative value of the fourth variable, as illustrated by
\begin{eqnarray*}
   && R_J(2,3,4,-0.5) = 0.24723\;81970\;3052\\
   && R_J(2,3,4,-5) = -0.12711\;23004\;2964
\end{eqnarray*}
In this example \(R_J\) vanishes when the fourth variable is approximately
\(-1.2552\), and cancellation between terms on the right side of (\ref{2.23a})
will lead to loss of significant figures.
\vspace{.1in}
\begin{eqnarray*}
   &&R_D(0,2,1) = 1.7972\;10352\;1034\\
   &&R_D(2,3,4) = 0.16510\;52729\;4261\\
   &&R_D(i,-i,2) = 0.65933\;85415\;4220\\
   &&R_D(0,i,-i) = 1.2708\;19627\;1910\;+\;i\;2.7811\;12015\;9521\\
   &&R_D(0,i-1,i) = -1.8577\;23543\;9239\;-\;i\;0.96193\;45088\;8839\\
   &&R_D(-2-i,-i,-1+i) = 1.8249\;02739\;3704\;-\;i\;1.2218\;47578\;4827\\
\\
   &&R_G(0,16,16) = 2E(0) = \pi = 3.1415\;92653\;5898\\
   &&R_G(2,3,4) = 1.7255\;03028\;0692\\
   &&R_G(0,i,-i) = 0.42360\;65423\;9699\\
   &&R_G(i-1,i,0) = 0.44660\;59167\;7018\;+\;i\;0.70768\;35235\;7515\\
   &&R_G(-i,i-1,i) = 0.36023\;39218\;4473\;+\;i\;0.40348\;62340\;1722\\
   &&R_G(0,0.0796,4) = E(0.99) = 1.0284\;75809\;0288
\end{eqnarray*}
\indent
Consistency checks that do not use external information are furnished by the 
following equations, in which \(x,y,p\) are positive, \(\l \m = xy\), and 
\( |{\rm ph}\,\l | < \pi\)\,:
\begin{equation}\label{3.1}
   R_F(x+\l,y+\l,\l) + R_F(x+\m,y+\m,\m) = R_F(x,y,0)\,;
\end{equation}
\begin{equation}\label{3.2}
   R_C(\l,x+\l) + R_C(\m,x+\m) = R_C(0,x)\,;
\end{equation}
\begin{equation}\label{3.3}
   R_J(x+\l,y+\l,\l,p+\l) + R_J(x+\m,y+\m,\m,p+\m) = R_J(x,y,0,p)
   - 3\,R_C(a,b)\,,
\end{equation}
where
\begin{equation}\label{3.4}
   a = p^2(\l+\m+x+y)\,,\quad b = p(p+\l)(p+\m)\,,\quad b-a = p(p-x)(p-y)\,;
\end{equation}
\begin{equation}\label{3.5}
   R_D(\l,x+\l,y+\l) + R_D(\m,x+\m,y+\m) = R_D(0,x,y) - \frac{3}{y\sqrt{x+y+
   \l + \m}}\,.
\end{equation}
These equations are special cases of addition theorems.  Another consistency
check is furnished by
\begin{equation}\label{3.6}
   R_D(x,y,z) + R_D(y,z,x) + R_D(z,x,y) = \frac{3}{\sqrt{x}\,\sqrt{y}\,
   \sqrt{z}}\,,
\end{equation}
where \(x,y,z\) lie in the cut plane and each square root has positive real 
part.

\section{Other integrals}
Legendre's complete elliptic integrals \(K\) and \(E\) are given by
\begin{eqnarray}
K(k) & = & R_F(0,1-k^2,1), \label{4.1}\\
E(k) & = & 2R_G(0,1-k^2,1) \nonumber\\
& = & \frac{1-k^2}{3} [R_D(0,1-k^2,1) + R_D(0,1,1-k^2)], \label{4.2}\\
K(k) - E(k) & = & \frac{k^2}{3} R_D(0,1-k^2,1), \label{4.3}\\
E(k) - (1-k^2)K(k)  & = & \frac{k^2(1-k^2)}{3} R_D(0,1,1-k^2). \label{4.4}
\end{eqnarray}
In Legendre's incomplete integrals we shall use the abbreviation 
\(c = {\rm csc}^2\phi = 1/\!\sin^2\!\phi \)\,:
\begin{eqnarray}
   F(\phi ,k) &=& (\sin \phi) R_F(\cos^2\! \phi,1-k^2 \sin^2\! \phi,1)
     = R_F(c-1,c-k^2,c)\,, \label{4.5}\\
   E(\phi ,k) &=& R_F(c-1,c-k^2,c) - \frac{k^2}{3} R_D(c-1,c-k^2,c)\,,
     \label{4.6}\\
   \Pi (\phi ,k,n) &=& \int_0^{\phi} (1+n\sin^2\theta)^{-1} 
     (1-k^2\sin^2\theta)^{-1/2}\,d\theta \nonumber\\
     &=& R_F(c-1,c-k^2,c) - \frac{n}{3}\,R_J(c-1,c-k^2,c,c+n)\,. \label{4.7}
\end{eqnarray}
Some related integrals are
\begin{eqnarray}
   D(\phi ,k)&=&\int_0^{\phi} sin^2\theta\, (1-k^2 sin^2\theta)^{-1/2}\,
     d\theta\ = \frac{1}{3}\,R_D(c-1,c-k^2,c)\,, \label{4.8}\\
   K(k)Z(\b ,k)&=&\frac{k^2}{3}\,\sin \b \cos \b\, \sqrt{1-k^2 \sin^2\!\b}\,
     R_J(0,1-k^2,1,1-k^2 \sin^2\!\b)\,, \label{4.9}\\
   \Lambda_0(\b ,k)&=&\frac{2}{\pi}\, \frac{(1-k^2)\sin \b \cos \b}{\d}
   \nonumber\\
   &\cdot& \left[ R_F(0,1-k^2,1)+\frac{k^2}{3\d^2}\,R_J\left(0,1-k^2,1,
   1-\frac{k^2}{\d^2}\right) \right]\,, \label{4.10}
\end{eqnarray}
where \(\d = \sqrt{1-(1-k^2)\sin^2\!\b}\).  The function \(Z(\b ,k)\) is 
Jacobi's zeta function \cite[140.03]{BF71}, and \(\Lambda_0(\b ,k)\) is 
Heuman's lambda function \cite[150.01]{BF71}.

Bulirsch's elliptic integrals \cite{Bu65} are 
\begin{eqnarray}
   el\,1(x,k_c) &=& x\,R_F(1,1+k_c^2 x^2,1+x^2)\,, \label{4.11}\\
   el\,2(x,k_c,a,b)&=& ax\,R_F(1,1+k_c^2 x^2,1+x^2) \nonumber \\
     &+& \frac{1}{3}(b-a)x^3\,R_D(1,1+k_c^2 x^2,1+x^2)\,, \label{4.12}\\
   el\,3(x,k_c,p)&=& x\,R_F(1,1+k_c^2 x^2,1+x^2) \nonumber\\
     &+& \frac{1}{3}(1-p)x^3\,R_J(1,1+k_c^2 x^2,1+x^2,1+px^2)\,,\label{4.13} \\
   cel(k_c,p,a,b)&=& a\,R_F(0,k_c^2,1) + \frac{1}{3}(b-pa)R_J(0,k_c^2,1,p)\,.
   \label{4.14}
\end{eqnarray}
\indent In the real domain, inverse circular and inverse hyperbolic 
functions are expressed in terms of \(R_C\) by
\begin{eqnarray}
      \ln \left(\frac{x}{y}\right) & = &  (x-y)R_{C}\left(
                  \left({\frac{x+y}{2}} \right)
                  ^{2},xy \right),\qquad  x > 0, \label{4.16}\\
    \arctan (x/y) & = & xR_{C}(y^2,y^2+x^2),\qquad  -\infty < x < \infty,
    \label{4.17}\\
    \hbox{arctanh}(x/y) &= &xR_{C}(y^2,y^2-x^2),\qquad  -y < x < y,
    \label{4.18}\\
    \arcsin(x/y) & = &xR_{C}(y^2-x^2,y^2),\qquad -y \leq x \leq y,
    \label{4.19}\\
    \hbox{arcsinh}(x/y) &= &xR_{C}(y^2+x^2,y^2),\qquad -\infty < x < \infty,
    \label{4.20}\\
    \arccos(x/y) & = &(y^2-x^2)^{\frac{1}{2}}R_{C}(x^{2},y^2),\qquad
    0 \leq x \leq y, \label{4.21}\\
    \hbox{arccosh}(x/y) &= &(x^2-y^2)^{\frac{1}{2}}R_{C}(x^2,y^2),\qquad
    x \geq y, \label{4.22}
\end{eqnarray}
where \(y > 0\) in each case.
These equations remain valid in the complex domain provided that the variables 
of \(R_C\) satisfy the conditions accompanying (\ref{1.3}).  If \(y=1\) the 
function multiplying \(R_C\) shows in each case the asymptotic behavior as the 
left side tends to 0.

Many elliptic integrals of the form 
\begin{equation}\label{4.15}
   \int_y^x \prod_{i=1}^n (a_i + b_i t)^{p_i/2} dt\,,
\end{equation}
where \(p_1,\ldots,p_n\) are integers and the integrand is real, are 
reduced in \cite{Ca87}-\cite{Ca92} to the integrals in Section 2.  Use of the 
algorithm for \(R_F\) is illustrated by numerical computation of the integral
\begin{equation}\label{4.23}
   I = \int_y^x \frac{dt}{\sqrt{(f_1 + 2g_1 t + h_1 t^2)(f_2 + 2g_2 t + 
   h_2 t^2)}}\,,
\end{equation}
where all quantities are real, \(x > y\), the two quadratic (or, if \(h_i=0\),
linear) polynomials are positive on the open interval of integration, and 
their product has at most simple zeros on the closed interval.  Let
\begin{eqnarray}
   q_i(t)&=&f_i + 2g_i t + h_i t^2,\qquad i = 1,2\,,\label{4.24}\\
   (x-y)U&=& \sqrt{q_1(x)q_2(y)} + \sqrt{q_1(y)q_2(x)}\,,\label{4.25}\\
   T&=& 2g_1 g_2 - f_1 h_2 - f_2 h_1\,,\label{4.26}\\
   V&=& 2\sqrt{(g_1^2 - f_1 h_1)(g_2^2 - f_2 h_2)}\,. \label{4.27}
\end{eqnarray}
Then
\begin{equation}\label{4.28}
   I = 2\,R_F(U^2 + T + V,\,U^2 + T - V,\,U^2)\,.
\end{equation}
Except for notation this is the same as \cite[(34)]{Ca77}.  If the interval of 
integration is infinite, U is obtained by taking a limit; for example, if 
\(x = +\infty\)\,, then \(U = \sqrt{h_1 q_2(y)} + \sqrt{q_1(y)h_2}\)\,.
If exactly one of
\(q_1\) and \(q_2\) has conjugate complex zeros, then \(V\) is pure imaginary;
otherwise the variables of \(R_F\) are real.
The algorithm for \(R_F\) in this paper can be used in both cases.
However, if both polynomials have conjugate complex zeros, the quadrilateral 
with the zeros as vertices has diagonals that must not intersect at an interior
point of the interval of integration.  If they do, the integral must be split
into two parts at the point of intersection.  This restriction is discussed in 
\cite[\S 4]{Ca77} and removed by a Landen transformation in \cite{Ca92}.

\end{document}